\documentclass[10pt]{amsart}
\usepackage{xspace}
\usepackage{amssymb}

\usepackage[hyperindex, bookmarksnumbered,
 plainpages]{hyperref}

\newcommand\dcl\DeclareMathOperator
\dcl\Sp{Specm} \dcl\depth{depth} \dcl\im{Im}
\dcl\codim{codim} \dcl\Y{Y} \dcl\gl{gl}
\dcl\U{U} \dcl\T{T} \dcl\qdet{qdet}
\dcl\sgn{sgn} \dcl\gr{gr} \dcl\diag{diag}
\dcl\g{g} \dcl\C{\mathbb C} \dcl\dd{{\mathrm
d}}

\newcommand\CM{\mbox{\textrm{CM}}\xspace}
\newcommand\CI{\mbox{\textrm{CI}}\xspace}

\newcommand\sC{{\mathsf C}}
\newcommand\sU{{\mathsf U}}
\newcommand\sV{{\mathsf V}}
\newcommand\sX{{\mathsf X}}

\newcommand\K{{\mathbf K}}
 
\newcommand\sg{{\mathsf g}}
\newcommand\sh{{\mathsf h}}
\newcommand\osg{\overline{\mathsf g}}

\renewcommand\sf{{\mathsf f}}

\newcommand\og{\overline{g}}
\newcommand\oh{\overline{h}}
\newcommand\ou{\overline{u}}

\newcommand\Ga{{\Gamma}}
\newcommand\La{{\Lambda}}
\newcommand\oGa{\overline{\Gamma}}
\newcommand\oLa{\overline{\Lambda}}
\newcommand\oU{\overline{U}}
\newcommand\oI{\overline{I}}

\renewcommand\k{{\Bbbk}}

\newcommand\bZ{{\mathbb Z}}

\newcommand\leqs{\leqslant}
\newcommand\geqs{\geqslant}

\newcommand\ds{\displaystyle}

\newtheorem{theorem}{Theorem}
\newtheorem{lemma}{Lemma}[section]
\newtheorem{corollary}{Corollary}
\newtheorem{proposition}{Proposition}

\title{KOSTANT THEOREM FOR SPECIAL FILTERED ALGEBRAS}

\author{Vyacheslav Futorny}

\address{\noindent School of Mathematics and Statistics, University of Sydney \\
NSW 2006, Sydney, Australia}
\email{futorny@maths.usyd.edu.au, futorny@ime.usp.br}

\author{Serge Ovsienko}

\address{\noindent
Faculty of Mechanics and Mathematics\\
Kiev Taras Shevchenko University\\
Vla\-di\-mir\-skaya 64, 01033, Kiev, Ukraine}
\email{ovsienko@sita.kiev.ua }

\begin{document}

\begin{abstract} A famous result of Kostant  states that the universal
enveloping algebra of a semisimple complex Lie
algebra is a free module over its center. We
prove an analogue of this result for the class
of  special filtered algebras and
apply it to show the freeness over its center of
the restricted Yangian and the universal
enveloping algebra of the restricted current algebra,
associated with the general linear Lie algebra.
\end{abstract}

\maketitle

\vspace{0,2cm}

\subjclass{Mathematics Subject Classification
13A02, 16W70, 17B37, 81R10}

\section {Introduction}

A well-known theorem of Kostant \cite{K} says
that for a complex semisimple Lie algebra $\g$
the universal enveloping algebra $\U(\g)$ is a
free module over its center.  For the Yangian
$\Y(\gl_n)$ of the general linear Lie algebra
$\gl_n$ the freeness over the center  was shown
in \cite{MNO}.
 Another example of such
situation gives the Gelfand-Tsetlin subalgebra
$\Ga$ of $\U(\gl_n)$ (see \cite{DFO},
\cite{O1}). It
 was shown in \cite{O2} that $\U(\gl_n)$ is a free left (right)
 $\Ga-$module. The knowledge of the freeness of a given algebra over its certain
subalgebra often allows to establish other
important results, e.g. to find the
annihilators of Verma modules \cite{Di}, to
prove the finiteness of the number of liftings
to an irreducible module over $\U(\gl_n)$ from
a given character of the Gelfand-Tsetlin
subalgebra (\cite{O2}), etc.

In  \cite{O2} was introduced a technique which
generalizes Kostant methods (\cite{K}, see also
\cite{G1}) and which allows to study the
universal enveloping algebras of Lie algebras
as modules over its certain commutative
subalgebras. In the present paper we develop a
graded version of the technique from \cite{O2}
which can be applied  to a large class
of filtered associative  algebras. In particular, we apply our
result to the restricted Yangians $\Y_p(\gl_n)$
of any level $p$ and to the universal
enveloping algebra of restricted current
algebras  of any level $m$, associated with the
general linear Lie algebra, and show that these
algebras are free as modules over their
centers. A particular case of this statement for
$\Y_2(\gl_n)$  was considered by
Geoffriau (\cite{G1}) who  showed that the
universal enveloping algebra of a Takiff
algebra is free over its center.

The structure of the paper is the following. In
Section \ref{section-regular-sequences} we
collect necessary  known facts about Koszul
complexes and complete intersections. In
Section \ref{section-filtered-algebras} we
consider a filtered algebra $U$ containing some
commutative subalgebra $\Ga$. We find
sufficient conditions that guarantee the
lifting of every character of $\Ga$ to a simple
$U$-module. In Section
\ref{section-an-analogoue-of-Kostant-theorem}
we establish the key result for special
filtered algebras (Theorem
\ref{theorem-main-analogue-kostant-theorem})
showing their freeness over  commutative
subalgebras generated by a sequence of elements
whose graded images form a complete
intersection for the associated graded algebra
(equivalently, the corresponding characteristic
variety is equidimensional of minimal possible
dimension). This generalizes the Kostant
theorem for semisimple Lie algebras. In Section
\ref{section-applications} we apply Theorem
\ref{theorem-main-analogue-kostant-theorem} to
the restricted Yangian of any level $p$ for
$\gl_n$ and to the universal enveloping
algebras of the restricted current algebras of
any level $m$ for $\gl_n$  (Theorem
\ref{theorem-application-to-restricted-yangians}
and
\ref{theorem-application-to-current-algebras}).

Throughout the paper we fix an algebraically
closed field $\k$ of characteristic $0$. For an
affine $\k$-algebra $\La$ (i.e. associative and
commutative finitely generated ${\k}$-algebra)
we denote by $\Sp \La$ the variety of all
maximal ideals of $\La$. If $\La$ is a
polynomial algebra in $n$ variables then we
identify   $\Sp \La$ and $\k^n$. $\dim \La$
($=\dim \Sp \La$) means the Krull dimension.
For an ideal $I\subset {\La}$ denote by
$\sV(I)\subset$ $\Sp{\La}$ the set of all
zeroes of $I$, $\sV(I)=\{\mu\in \Sp \La|
I\subset\mu\}$. If $I$ is generated by $f_1,
\ldots, f_s$ then we write $I=(f_1, \ldots,
f_s)$ and $\sV(I)=\sV(f_1, \ldots, f_s)$.  The
word "graded" always means "positively graded".

\section {Koszul complex and complete intersections}
\label{section-regular-sequences}

The main sources of references related this
chapter are \cite{Ei}, \cite{Ma} and \cite{B}.

Let  $U$ be an associative (not necessary
commutative) algebra and $M$  an $U-\k[X_1,$
$\dots,$ $X_t]$ bimodule. Recall (\cite{B}, \S
9) that  the associated \emph{Koszul complex}
$\K_.\ (=\K_.(M))$ of left $U$-modules
 is defined  as follows. Let
$e_1,\dots,e_t$ be a standard basis of  $\k^t$.
Set $\K_i=0$ for $i<0$ and $\K_i=\ds
M\otimes_\k\bigwedge^i_\k\k^t$ for each
$i\geqs0$.  Define the differential
$d_i:\K_i\to \K_{i-1}$ on the element
 $m\otimes(e_{j_1}\wedge \dots \wedge e_{j_i})$ of $\K_i$, $1\leqs
j_1 <\dots< j_i\leqs t$, $m\in M$ as:

\begin{equation}\label{equ-formula-in-Koshul}d_i(m\otimes (e_{j_1}\wedge
 \dots \wedge e_{j_i}))=\sum_{k=1}^{i} (-1)^{k-1} m\cdot X_{j_k}\otimes
(e_{j_1}\wedge \dots \widehat{e}_{j_k} \dots\wedge e_{j_i}),
\end{equation}
where $\widehat{e}_{j_k}$ means the omission of
$e_{j_k}$.  Let $\sg=\{g_1,\ldots$,
$g_t\}\subset U$ be a sequence of mutually
commuting elements, $M=U$ with the natural
structure of left $U$-module and $a\cdot X_i=a
g_i$, $a\in U$, $i=1,\dots,t$. We denote the
corresponding Koszul complex by $\K_.(\sg,U)$.

The sequence $\sg$ is called a \emph{complete
intersection for $U$} (we say simply
\emph{complete intersection} if $U$ is fixed)
provided for  $\K_.(\sg,U)$ holds
$H_i(\sg,U)=0$ for $i\ne 0$ and
$H_0(\sg,U)\ne0$, where $H_i(\sg,U)$ is the
$i$-th homology of  $\K_.(\sg,U)$. Note that
$H_0(\sg,U)\simeq$ $ U/(Ug_1+\dots +Ug_t)$.

The sequence $g_1,\dots,g_t$ in an affine
algebra $\La$ is called \emph{regular},
provided the class of $g_i$ is not a zero
divisor and not invertible in
$\La/(g_1,\dots,g_{i-1})$ for any
$i=1,\dots,t$. A regular sequence is a complete
intersection for $\La$.

An affine algebra $\La$ is called a
\emph{complete intersection, shortly \CI},
provided $\La\simeq$ $ A/I$, where $A$ is a
polynomial algebra ${\k}[X_1,\ldots,X_n]$ and
the ideal $I$ is generated by a regular
sequence $\sg=\{g_1,\dots,g_t\}$ in $A$.

An affine algebra $\La$ is called
\emph{Cohen-Macauley, shortly \CM,} provided
for every ideal $I\subset \La$ holds $\depth I
=\codim I$, where $ \depth I$ is a length of
maximal regular sequence, contained in $I$ and
$\codim I$ is a minimal codimension of
irreducible components of $\sV(I)$ in $\Sp
\La$. If $g_1,\dots,g_t$ is a sequence of
elements in a \CM algebra $\La$ such that  $\sV(I)$ is equidimensional
 for
$I=(g_1,\dots,g_t)$ and $\codim I=t$ then $\La/I$ is \CM
(\cite{Ei}, Proposition 18.13). In particular,
any \CI algebra is a \CM algebra. \emph{We will
assume that a \CI algebra $\La\simeq A/I$ is
reduced, i.e. $I$ is a radical ideal.}

\begin{proposition}\label{proposition-bh}
Let $\La$ be an affine algebra of dimension $n$
and $\sg=\{g_1,\dots,g_t\}$ a sequence of
elements of $\La$, $0\leqs t\leqs n$.

\begin{enumerate}
\item(\cite{B}, \S 9, Proposition 
6)\label{item-proposition-linear-change-keeps-complete-intersection}
Let $\sg$ be a complete intersection for $\La$,
$L\in GL_t(\La)$ and $(h_1, $ $\ldots,$ $
h_t)=$ $ (g_1, \ldots, g_t)\cdot L$. Then the
sequence $h_1, \ldots, h_t$ is a complete
intersection for $\La$.

\item(\cite{Ma},
Theorem 16.5,(ii), \cite{B}, \S 7, Corollary
2) \label{item-proposition-regular-homogeneous-equiv-complete-intersection}
If $\La$ is graded and $g_1,\dots,g_t$ are
homogeneous then $\sg$ is regular if and only
if $\sg$ is a complete intersection.

\item(\cite{Ma1},
Theorem 28) \label{item-subsequeqnce-is-complete-intersection}
Let $\La$ be a graded  algebra and $\sg$ be
a complete intersection for $\La$ consisting
of homogeneous elements. Then any subsequence
of $\sg$ is a complete intersection.

\item\label{item-proposition-equidimensional-hohmogeneous_is-regular}
Let $\La$ be a graded \CM algebra, $\sg$
consists of homogeneous elements and $\codim
I=t$, $I=(g_1,\dots,g_t)$. Then the sequence
$\sg$ is regular.

\item\label{item-proposition-equidimensionality-is-equivalent-to-complete-intersection}
Let $\La$ be a \CM algebra. Then  $\sg$ is a
complete intersection for $\La$ if and only if
the variety $\sV(g_1, \ldots, g_t)$ is
equidimensional of dimension $n-t$.

\item\label{item-strong-equidimensionality-is-equivalent-to-regularity}
Let $\La$ be a \CM algebra. Then  $\sg$ is a
regular sequence if and only if every variety
$\sV(g_1, \ldots, g_i)$ is equidimensional of
dimension $n-i$,  $i=1,\dots,t$.

\end{enumerate}
\end{proposition}

\begin{proof} It follows from
 \cite{Ma}, Theorem 16.8, that for a sequence $\sg=\{g_1, \ldots, g_t\}$
 in $\La$, such that
$I=(g_1,\dots, g_t)\ne\La$ holds $\depth I
=t-q$, where  $\ds
q=\max_i\{H_i(\g,\La)\ne0\}$. Hence $\sg$ is a
complete intersection if and only if $\depth
I=t$ and, in assumption $\La$ is \CM, if and
only if $\codim I=t$. On the other hand $\codim
I\geqs t$ and we have the equality if and only
if $\sV(I)$ is equidimensional of minimal
possible dimension $n-t$.
 This implies immediately
(\ref{item-proposition-equidimensional-hohmogeneous_is-regular})
and
(\ref{item-proposition-equidimensionality-is-equivalent-to-complete-intersection}).
In
(\ref{item-strong-equidimensionality-is-equivalent-to-regularity})
the statement "only if"  is obvious. The
statement "if" we prove by induction on $i$.
The base of induction $i=0$ is obvious. Denote
$\La_i=\La/(g_1,\dots,g_{i})$ and let $g$ be
the class of $g_{i+1}$ in $\La_i$. Obviously
$g$ is not invertible. Let $J=\La_i g$. Then in
$\La_i$ holds $\codim J=\depth J=1$ and hence
$J$ contains an element which is not a zero
divisor implying that $g$ is not a zero
divisor.
\end{proof}

We have the following  immediate corollary of
Proposition \ref{proposition-bh},
(\ref{item-proposition-equidimensionality-is-equivalent-to-complete-intersection}).

\begin{lemma}\label{lemma-substituion-0-in-regular-gives-regular}
Let  ${A}=$ ${\k}[X_1,\ldots,X_n]$ be the
polynomial algebra, $G_1, \ldots, G_t\in A$. A
sequence $X_1,$ $ \ldots,$ $ X_r,$ $ G_1,$ $
\ldots,$ $ G_t$
  is a complete intersection for $A$ if and only
if the sequence $g_1, \ldots, g_t$ is a
complete intersection for $\k[X_{r+1}, \ldots,
X_n]$, where $g_i(X_{r+1},$ $ \ldots,$ $ X_n)=$
$G_i(0,$ $ \ldots,$ $ 0,$ $ X_{r+1},$ $
\ldots,$ $ X_n)$, $i=1,\dots,t$.
\end{lemma}

\section{Filtered  algebras}
\label{section-filtered-algebras}

Let $U$ be a \emph{filtered} algebra, i.e an
associative algebra over ${\k}$, endowed with
an increasing exhausting filtration $\{U_{i}
\}_{i\geqs0}$, $U_{0} =$ ${\k}$, $U_i
U_j\subset U_{i+j}$, $U=\cup_{i\geqs 0} U_i$.
For $u\in U_i\setminus U_{i-1}$ set $\deg u=i$
($U_{-1} =$ $\{ 0 \}$). Let $\oU=\gr U$ be the
associated graded algebra ${\oU}=$ $\ds
\bigoplus_{i=0}^{\infty} U_{i }/U_{i-1 }$. For
$u\in U$ denote by $\ou$ its image in $\oU$ and
for a subset $S\subset U$ denote $S_i=S\cap
U_i$, ${\overline  S}=$ $\{ {\overline
s}\,|\,s\in S\}$ $\subset {\oU}$. Set
$\oU_{(i)}=U_i/U_{i-1}$ and for any $T\subset
\oU$ denote $T_{(i)}=T\cap \oU_{(i)}$. Given a
graded algebra $U=\bigoplus_{i\geqs0} U_{(i)}$,
$U_{(0)}=\k$ we always assume in $U$ the
associated filtration $\{U_i= \bigoplus_{j\leqs
i} U_{(j)}, i\geqs 0\}$ identifying $U$ and
$\oU$.

\begin{lemma}\label{lemma-keeping-the-regularity}
Let $\sg=$ $\{g_1,$ $\dots,$ $g_t\}$ be a
sequence of mutually commuting elements of a
filtered algebra $ U$ such that
$\osg=\{{\og}_1,\dots,{\og}_t\}$ is a complete
intersection for $\oU$. Then

\begin{enumerate} \item\label{item-lemma-keeping-the-complete-intersection} $\sg$ is a complete
intersection for $U$;

\item\label{item-lemma-keeping-graduation-of-ideal-coincides-with-generated-by-graduations}
if $I=U g_1+\dots+Ug_t$ then
$\oI=\oU{\og}_1+\dots+\oU{\og}_t$.
\end{enumerate}
\end{lemma}

\begin{proof}
Endow  the complex $\K_.=\K_.(\sg,U)$ with a
filtration $F$:
$$\dots \subset F^{p-1}\K_i\subset F^{p}\K_i\subset
F^{p+1}\K_i\subset\dots\subset \K_i,\text{
where } p,i\in\bZ,\text{ by setting }$$
 $$ u \otimes (e_{j_1}\wedge  \dots
\wedge e_{j_i})\in F^p\K_i \text{ if and only
if } \deg u + \deg g_{j_1}+\dots+\deg
g_{j_i}\leq p.
$$ $u\in U$, $1\leqs j_1<\dots< j_i\leqs t$.
This filtration is exhausting, bounded below
and turns $\K_.$ into a filtered complex.
Analogous grading turns $\K_.(\osg,\oU)$ into a
complex of graded modules. By the conditions
$H_i(\osg,\oU)=0$ for $i\ne 0$ and the $p$-th
graded component of $0$-homology is $\ds
H_0(\osg,\oU)_{(p)}=\oU_{(p)}/\sum_{i=1}^t\oU_{(p-\deg
g_i)}\og_i.$

Consider associated with the filtration $F$
spectral sequence $(E^r,$ $d^r)$ with the first
term  $(E^1,d^1)$,
$E^1_{p,q}=H_{p+q}(F^{p+1}\K_./F^{p}\K_.)$
(\cite{W}, 5.4). The  complex of graded
$U$-modules $\{F^{p+1}\K_./F^{p}\K_.\}_{p\in
\bZ}$ is isomorphic to $\K_.(\osg,\oU)$, hence,
the nonzero components of $E^1$ are among
$E^1_{p,-p},p\in\bZ$ (they are the graded
components of $H_0(\osg,\oU)$). Besides all
differentials in $E^1$ equal $0$. It gives us
$E^1\simeq E^2\simeq \dots\simeq E^\infty$.
Since $(E^r,d^r)$ converges to $H(\sg,U)$ it
follows that all homologies $H_i(\sg,U)=0$,
$i\ne 0$ which proves the statement
(\ref{item-lemma-keeping-the-complete-intersection}).
 To prove the  statement
(\ref{item-lemma-keeping-graduation-of-ideal-coincides-with-generated-by-graduations})
we consider subcomplex $\K'_.(\sg,U)\subset
\K_.(\sg,U)$ with $\K'_i(\sg,U)=\K_i(\sg,U)$,
$i\ne 0$, $\K'_0(\og,U)=0$, and a subcomplex
$\K'_.(\osg,\oU)\subset \K_.(\osg,\oU)$ with
$\K'_i(\osg,\oU)=\K_i(\osg,\oU)$, $i\ne 0$,
$\K'_0(\osg,\oU)=0$. The only nonzero
homologies of these complexes are:
$H_1(\K'_.(\osg,\oU))\simeq U\og_ 1+\dots+
U\og_t$ and  $H_1(\K'_.(\sg,U))\simeq I$
 respectively. The  filtration on $\K_.$
induces a filtration on $\K'_.$. As above,
considering the corresponding spectral
sequences  we conclude that
$\oI=U{\og}_1+\dots+ U{\og}_t$.
\end{proof}

In commutative case we use a  more precise
statement.

\begin{lemma}\label{lemma-regularity-filtered-affine-from-regularity}
Let $\La$ be a filtered affine algebra
 ($\La=\oLa$) and
$\sg=\{g_1, \ldots, g_t\}\subset \La$  a
sequence such that $\osg$  is a complete
intersection for $\oLa$. Then $\sg$ is a
regular sequence.
\end{lemma}

\begin{proof}
By Proposition \ref{proposition-bh},
(\ref{item-proposition-regular-homogeneous-equiv-complete-intersection})
and
(\ref{item-subsequeqnce-is-complete-intersection})
the sequence $\osg$ and any its subsequence are
regular. Following Lemma
\ref{lemma-keeping-the-regularity},
(\ref{item-lemma-keeping-graduation-of-ideal-coincides-with-generated-by-graduations}),
for the ideal $I_k=(g_1,\dots,g_k)$ holds
$\oI_k=(\og_1,\dots,\og_k)$ for any
$k=1,\dots,t$.  Let $g\not\in (g_1,\dots,g_i)$
be an element of minimal degree, such that
$g_{i+1}g\in (g_1,\dots,g_i)$. Then $\og_{i+1}
\og\in (\og_1,\dots,\og_i)$. By regularity of
$\osg$ we obtain that $\og\in (\og_1,\dots,\og_i)$ and
hence there exists $h\in(g_1,\dots,g_i)$  such
that $\deg (g-h) < \deg g$ and
$g_{i+1}(g-h)\in(g_1,\dots,g_i)$. The obtained contradiction proves
the lemma.
\end{proof}

\begin{lemma}{}\label{lemma-regular-in-complete-intersection}
Let  $\La$ be a graded affine algebra,
$\{\La_i,i\in \bZ\}$  the associated
filtration, $\sg=$ $\{g_1, $ $\ldots,$ $ g_t\}$
a complete intersection for $\La$ consisting of
homogeneous  elements, $\Ga=\k[g_1, \ldots,
g_t]$,  $\mu=$ $(\mu_1,$ $ \ldots, $
$\mu_t)\in$ $\La^t$ such that $\deg \mu_i<$
$\deg g_i$, $i=1,\dots,t$,
$I_{\mu}=(g_1-\mu_1,\ldots, g_t-\mu_t)$. Then

 \begin{enumerate}

\item \label{item-poly-min-present}
$I_{\mu}\cap \La_m=$ $\ds {\sum_{i=1}^t}{\La}_{
m-d_i}(g_i- \mu_i)$ where $d_i=$ $\deg g_i$,
$i=1,\ldots,t$.

\item\label{item-poly-shifted-is-regular} The
sequence $\{g_1-\mu_1,\ldots,g_t-\mu_t\}$ is a
complete intersection for ${\La}$, moreover, it
is regular.

\item\label{item-poly-graded-does-not-depend}
$\oI_{\mu}$ is an ideal generated by
$g_1,\ldots, g_t$, in particular $I_{\mu}\ne
{\La}$.

 \item\label{item-poly-map-is-epi} If
$\La$ is \CM, then the regular map $p_{\sg}:$
$\Sp{\La}{\longrightarrow}$ $\Sp{\Ga}=$ $\k^t$
induced by the inclusion
$i_{\sg}:{\Ga}\hookrightarrow{\La}$ is an
epimorphism and $\dim p_{\sg}^{-1}(\mu)=$ $n-t$
for any $\mu\in$ $\Sp{\Ga}$.
\end{enumerate}\end{lemma}

\begin{proof} Statements (\ref{item-poly-min-present}) and
 (\ref{item-poly-graded-does-not-depend}) follow from Lemma
\ref{lemma-keeping-the-regularity},
(\ref{item-lemma-keeping-graduation-of-ideal-coincides-with-generated-by-graduations}),
while  the statement
(\ref{item-poly-shifted-is-regular}) follows
from Lemma
\ref{lemma-regularity-filtered-affine-from-regularity}.
Applying to  $\mu=(\mu_1,\dots,\mu_t)\in\k^t$
(\ref{item-poly-shifted-is-regular}) of this
lemma and Proposition \ref{proposition-bh},
(\ref{item-proposition-equidimensionality-is-equivalent-to-complete-intersection})
we obtain, that
$p_{\sg}^{-1}(\mu)=\sV(g_1-\mu_1, \ldots,
 g_t-\mu_t)$  is equidimensional of dimension
$n-t$. This implies
(\ref{item-poly-map-is-epi}).

\end{proof}

\begin{proposition}{}\label{proposition-lem_2.5}
Let $U$ be a filtered algebra,
$\sg=\{g_1,\dots,g_t\}$ be a sequence of
mutually commuting elements from $U$ such that
$\{{\og}_1$, $\ldots$, ${\og}_t\}$ is a
complete intersection for ${\oU}$ and $\mu_i\in
U$ be such, that $\deg g_i>\deg \mu_i$, $i=1,
\ldots, t$. Denote $\mu=(\mu_1, \ldots,
\mu_t)$,
 $I_{\mu}=$
$U(g_1-\mu_1)+\ldots$ $+U(g_t-\mu_t)$ and
assume that the elements
$\{g_i-\mu_i|i=1,\dots,t\}$ are mutually
commuting (e.g. it holds if $\mu_i\in \k$ for
all $i$). Then
\begin{enumerate}

\item\label{item_pbw_min_present} $I_{\mu}\cap
U_m=$ $\ds \sum_{i=1}^t U_{m-d_i}(g_i- \mu_i)$
where $d_i=$ $\deg g_i$.

\item\label{item_pbw_graded_does_not_depend}
${\oI}_{\mu}=$ ${\oU}{\og}_1 +\ldots +
{\oU}{\og}_t$, in particular $I_{\mu}\ne$ $U$.


\item\label{item_CI-and-smaller_exists_a_simple}
There exists a simple left $U-$module $M$
generated by $m\in M$ such that for every $g_i$
holds $g_i m=$ $\mu_i m$.

\item
\label{item_CI_for-character-exists_a_simple}For
any $\nu\in\Sp \Ga$ there exists a simple left
$U-$module $M$ generated by $m\in M$ such that
for every $\gamma\in\Ga$ holds $\gamma\cdot m=$
$\nu(\gamma) m$ (i.e. $\nu$ lifts to a simple
$U$-module).

\end{enumerate}\end{proposition}

\begin{proof}
Statements (\ref{item_pbw_min_present}),
(\ref{item_pbw_graded_does_not_depend}) follow
from  Lemma \ref{lemma-keeping-the-regularity},
(\ref{item-lemma-keeping-graduation-of-ideal-coincides-with-generated-by-graduations}).
To obtain from
(\ref{item_pbw_graded_does_not_depend}) the
statement
(\ref{item_CI-and-smaller_exists_a_simple}) we
consider a maximal  left ideal $\mathfrak m$ in
$U$, containing $I_{\mu}$, and set $M=$ $U /
\mathfrak m$, $m=1+\mathfrak m$. Statement
(\ref{item_CI_for-character-exists_a_simple})
follows from
(\ref{item_CI-and-smaller_exists_a_simple})
applied to $\mu_i=g_i(\nu)$, $i=1,\dots, t$.
\end{proof}

\begin{corollary}
\label{corollary-equidimensionality-implies-existence-of-representation}
Let $U$ be a  filtered algebra such that
the associated graded algebra $\oU$ is CM. If
$g_1,\dots,g_t$ are mutually commuting elements
of $\,U$ such that
$\sV(\og_1,\dots,\og_t)\subset\Sp \oU$ is
equidimensional of codimension $t$, then every
$\mu\in\Sp \k[g_1,\dots,g_t]$ lifts to a simple
$U$-module.
\end{corollary}

\begin{proof}
 Follows from Proposition
\ref{proposition-lem_2.5},
(\ref{item_CI_for-character-exists_a_simple})
and Proposition \ref{proposition-bh},
(\ref{item-proposition-equidimensionality-is-equivalent-to-complete-intersection}).
\end{proof}

\section{An analogue of Kostant theorem}
\label{section-an-analogoue-of-Kostant-theorem}

A filtered algebra $U$ such that $\oU$ is a
reduced \CI algebra will be called {\it special
filtered}. Due to the Poincar\'e-Birkhoff-Witt
theorem the universal enveloping algebra of a
finite-dimensional Lie algebra is special
filtered. In this section we prove the
following analogue of Kostant theorem
(\cite{K}) for special filtered algebras.

\begin{theorem}{}\label{theorem-main-analogue-kostant-theorem} Let $U$ be a
special filtered
 algebra,
 $g_1$, $\ldots$, $g_t\in U$
 mutually commuting elements such that
${\og}_1$, $\ldots$, ${\og}_t$ is a complete
intersection for ${\oU}$, ${\Ga}=$
${\k}[g_1,\ldots,g_t]$. Then $U$ is a free left
(right) ${\Ga}-$module.
\end{theorem}

We start from the following known fact

\begin{lemma}{}\label{lemma-min-intersection-is-open} Let $V$ be a finite-dimensional
space, $Y\subset V$ a subspace, $Gr_{r}(V)$ the
Grassmanian of $r$-dimensional subspaces of
$V$, $\La$ an affine algebra, $F:$
$\Sp{\La}{\longrightarrow}$ $Gr_{r}(V)$ a
regular map and $$d \ = \ \ds {\min_{\mu\in
\Sp{\La}}} \dim\,Y \cap F(\mu).$$ Then the set
$\sU=$ $\{\mu\in \Sp{\La}\,|$ $\dim Y\cap
F(\mu)=d \}$ is nonempty and open in $\Sp
{\La}$.
\end{lemma}

\begin{proof}
The function $\varphi:Gr_r(V)\to \bZ$,
$W\longmapsto \dim W\cap Y$ is upper
semi-continuous, i.e. for every $n\in \bZ$ the
set $G_n=\{W\in Gr_r(V)\,|\,\varphi(W)\leqs
n\}$ is open in $Gr_r(V)$. Then
$\sU=F^{-1}(G_d)$ is open and  the statement of
the lemma follows.
\end{proof}

Let $A=\k[X_1,\dots,X_n]$ be the polynomial
algebra. For $f\in A$ denote by $\dd( f)$ the
column
$\begin{pmatrix}\frac{\partial f}{\partial X_1}\\
\vdots\\ \frac{\partial f}{\partial
X_n}\end{pmatrix}\in A^n$ and for a sequence
$\sf=$ $\{f_1,\dots,f_k\}$ denote by $J(\sf)=
J(f_1,\dots,f_k)$ the $n\times k$
\emph{Jacobian  matrix} $\Big(\frac{\partial
f_j}{\partial X_i}\Big),$ $i=1,$ $\dots,$ $n;$
$ j=1,$ $\dots,$ $k$. If $x\in \k^n$ then
$r_{\sf}(x)$ will denote the rank of the matrix
$J(\sf)(x)$.
 If $\sf$ is a complete
intersection then the ideal
$I=(f_1,\dots,f_k)\subset A$ is radical if and
only if for every component $\sC$ of $\sV(I)$
there exists an open dense $\sU\subset \sC$ such that
  $r_{\sf}(x)=k$ for every $x\in \sU$
 (Theorem 18.15, \cite{Ei}).

\begin{lemma}
\label{lemma-complete-intersection-by-adding-allmost-any}
Let $\tilde{h}=\{h_1,\dots,h_s, g\}\subset A$, $s\geqs 0$, be a
complete intersection  such that
 $I=(h_1,$ $\dots,$ $h_s)$ is
 radical ideal. Consider the set $\sU$
consisting of $\mu\in\k$ such that the algebra
$A/I_\mu$  is reduced \CI, where $I_\mu=$
$(h_1,$ $\dots,$ $h_s,$ $g-\mu)$. Then $\sU$ is
open dense in $\k$.
\end{lemma}
\begin{proof}
Set  $\sX=\sV(I)$ and let
$\sX=\sC_1\cup\dots\cup\sC_N$ be the
decomposition of $\sX$ into  irreducible
components, $e:\sX\to \k$  the regular map,
$e(x)=g(x)$, $x\in \sX$ (hence
$\sV(I_\mu)=e^{-1}(\mu)$). For any $\sC_i$ the
image $e(\sC_i)\subset \k$ is either a point or
a nonempty open subset in $\k$. Since $\dim
e^{-1}(0) =\dim \sV(I_0)=n-(s+1)$ there exists
$\sC_i$ such that $e(\sC_i)$ is dense in $\k$.
If $\sU'\subset\k\setminus\ds \bigcup_{i,\dim
e(\sC_i)=0}e(\sC_i)$ consists of  points
$\mu$  such that $\dim e^{-1}(\mu)$ is maximal
possible,  then for $\mu\in\sU'$ the variety
$\sV(I_\mu)$ is equidimensional of dimension
$n-(s+1)$ and hence
 $h_1,\dots,h_s,g-\mu$ is a complete
intersection for $A$.

We can assume that the map $e|_{\sC_i}:\sC_i\to
\k$ is dominant for every $i=1,\dots,N$. Let
$\sh=\{h_1,\dots,h_s\}$. By the conditions for
any $i=1, \ldots, N$ there exists an open dense
set $\sU_{i}\subset \sC_i$ such that
$r_{\sh}(x)=s$ for any $x\in\sU_i$, i.e. there
exists a $s\times s$-minor $\Delta_i$ of
$J(\sh)$  which is nonzero on an open subset of
$\sC_i$.

Assume first that $\sX=\sC_1$ and
$\Delta=\Delta_1$ is invertible in $\La$. Let
$K$ be the field of fractions of $\La$. If
$r_{\tilde{\sh}}|_{\sX}\leqs s$ then the column
$\dd( g)|_{\sX}\in \La^n$ is a linear
combination of  $\dd( h_1)|_{\sX},\dots,\dd(
h_s)|_{\sX}\in\La^n$ with coefficients from
$K$. Moreover due to the invertibility of
$\Delta$ holds $\dd (g)|_{\sX}= \lambda_1\dd(
h_1)|_{\sX}+\dots+\lambda_s\dd (h_s)|_{\sX}$
for some $\lambda_1,\dots,\lambda_s\in\La$.
Lifting this equality in $A^n$ we obtain  that
$v= \ds \dd (g) -\sum_{i=1}^s f_i \dd( h_i)\in
I^n$ for $f_i\in A$, $f_i|_{\sX}=\lambda_i$.
Consider the element $h=\ds g-\sum_{i=1}^s f_i
h_i\in A$. Then $ \ds \dd(h)= v  -\sum_{i=1}^s
\dd(f_i)  h_i\in I^n$ and hence $\dd(
h)|_{\sX}=0$ implying that $h$, and therefore
$g$ is a constant on $\sX$. It follows that
$r_{\tilde{\sh}}(x)=s+1$ for an open subset in
$\sX$.

Now consider the general case and fix $i$,
$1\leqs i\leqs N$. Take a polynomial $f$ which
is nonzero on $\sC=\sC_i$ and vanishes on all
other components of $\sX$. Set
$\Delta=\Delta_i$. Consider a sequence
$\sh'=\{h_1,\dots,h_s, 1+f\Delta X_{n+1}\}$ in
$A'=\k[X_1,\dots,X_n,X_{n+1}]$ and the ideal
$I'=(h_1,$ $\dots,$ $h_s, $ $1+f\Delta
X_{n+1})$. Then the projection of
$\sC'=\sV(h_1,$ $\dots,$ $h_s,$ $ 1+f\Delta
X_{n+1})$ on the first $n$ coordinates
identifies $\sC'$ with an open subset
$\sC\setminus \sV(f\Delta)$, i.e. $\Delta$ is
invertible on $\sC'$. Proposition
\ref{proposition-bh},
(\ref{item-strong-equidimensionality-is-equivalent-to-regularity})
 implies  that $\sh'$ is a regular sequence. Moreover, the ring
$\La'=A'/I'$ is isomorphic to the localization
of $\La$ by $f\Delta|_{\sX}$, hence is reduced.
Then as above we can prove that the rank of
$$J(h_1,\dots, h_s, 1+f\Delta
X_{n+1}, g-\mu) =
\left(\begin{smallmatrix}{}\frac{\partial
h_1}{\partial X_1}&\dots&\frac{\partial
h_s}{\partial X_1}&\frac{\partial
(f\Delta)}{\partial X_1}X_{n+1} &
\frac{\partial g}{\partial X_1}\\
\vdots&\ddots&\vdots&\vdots&\vdots\\
\frac{\partial h_1}{\partial
X_n}&\dots&\frac{\partial h_s}{\partial
X_n}&\frac{\partial (f\Delta)}{\partial
X_n}X_{n+1}&
\frac{\partial g}{\partial X_n}\\
0&\dots&0&f\Delta&0\\
\end{smallmatrix}\right)$$ equals  $s+2$ on an open dense set $\sU\subset\sC'$.
Since $f\Delta$ is  invertible on $\sC'$ we
conclude that  the rank of $J(h_1,\dots,h_s,g)$
equals $s+1$ on $\sU$, hence the ideal
$I_{\mu}$ is radical. The statement follows.
\end{proof}

\begin{lemma}
\label{lemma-regular-over-ci-allmost-everywhere-reduced}

Let $\La$ be a reduced graded \CI algebra,
$\sg=\{g_1,\dots,g_t\}\subset\La$  a complete
intersection for $\La$ consisting of homogeneous elements,
$\Ga=\k[g_1, \dots,
g_t]$. Then there exists an open dense
${\sU}_{\sg}\subset$ $\Sp{\Ga}$ such that for
any $\mu\in$ ${\sU}_{\sg}$ the ideal
$I_{\mu}=(g_1-g_1(\mu),\dots,g_t-g_t(\mu))\subset$
$\La$ is radical.
\end{lemma}

\begin{proof}
Denote $\sX_{\sg}=\{\mu\in\Sp \Ga\,|\,I_\mu
\text{ is radical}\}$. The criterion of
radicality of $I_{\mu}$ before Lemma
\ref{lemma-complete-intersection-by-adding-allmost-any}
 implies that $\sX_{\sg}$  is a
constructable set. Hence
 it is enough to show that $\sX_{\sg}$ is
dense in $\Sp\Ga$. We proceed by induction on
$t$. The base of induction $t=1$ follows
immediately from Lemma
\ref{lemma-complete-intersection-by-adding-allmost-any}.
 Suppose that
 there exists an open dense
$\sU'\subset \k^{t-1}$ such that for any
$\mu'=(\mu_1,\dots,\mu_{t-1})\in\sU'$  the
ideal
$I_{\mu'}=(g_1-\mu_1,\dots,g_{t-1}-\mu_{t-1})$
is radical. By Lemma
\ref{lemma-regular-in-complete-intersection},
(\ref{item-poly-shifted-is-regular}) the
sequence $g_1-\mu_1,\dots,g_{t-1}-\mu_{t-1}$ is
regular and therefore the algebra
$\La_{\mu'}=\La/(g_1-\mu_1,\dots,g_{t-1}-\mu_{t-1})$
is a reduced CI. By Lemma
\ref{lemma-complete-intersection-by-adding-allmost-any}
for any such $\mu'$ there exists an open dense
$\sU_{\mu'}\subset \k$  such that for
$\mu_t\in\sU_{\mu'}$  the
 algebra
$\La_{\mu'}/(g_t-\mu_t)\simeq
\La/(g_1-\mu_1,\dots,g_{t-1}-\mu_{t-1},g_t-\mu_t)$
is a reduced \CI, therefore  the ideal
$(g_1-\mu_1,\dots, g_t-\mu_t)$ is radical.
Since
$\sU=\bigcup_{\mu'\in\sU'}\bigcup_{\mu_t\in\sU_{\mu'}}(\mu',\mu_t)\subset
\sX_{\sg}$ is dense in $\Sp \Ga$  the proof
follows.
\end{proof}

Till the end of this section we follow closely
the techniques from \cite{Di}, 8.2.3 and
\cite{O2}. In particular Lemmas
\ref{lemma-equality-modulo-ideal-and-equality-as-functions}
and \ref{lemma-freedom-in-commutative-case}
below are analogous to \cite{Di}, 8.2.1, 8.2.2.
Let $\La$ be an affine algebra. Following
\cite{Di} we say that $h_1,\ldots,h_k\in$
${\La}$ are \emph{linearly independent over an
ideal} $I\subset {\La}$, provided $h_1+I$,
$\ldots$, $h_k+I\in$ ${\La}/I$ are linearly
independent over ${\k}$.

\begin{lemma}{}\label{lemma-equality-modulo-ideal-and-equality-as-functions}
Let $\sg=\{g_1, \ldots, g_t\}$ be  a complete
intersection of homogeneous elements for a
reduced graded \CI algebra  $\La$, ${\Ga}=$
${\k}[g_1,\ldots,g_t]$ and let
$h_1,\ldots,h_k\in$ ${\La}$ be linearly
independent over $I=$ $(g_1,\ldots,g_t)$. Then
\begin{enumerate}
\item\label{item-linearly-indep-over-i-mu}
There exists an open dense set ${\sU}_1\subset$
$\Sp{\Ga}$ such that for every $\mu\in {\sU}_1$
 $h_1,$ $\ldots,$ $h_k$ are
linearly independent over
$I_{\mu}=(g_1-g_1(\mu),\dots,g_t-g_t(\mu))$.

\item\label{item-linearly-indep-as-functions}
There exists an open dense ${\sU}_2\subset$
$\Sp{\Ga}$ such that the restrictions of $h_1,
\ldots, h_k$ on $\sV(I_{\mu})$ are linearly
independent over ${\k}$ for each $\mu\in
{\sU}_2$.

\end{enumerate}\end{lemma}

\begin{proof}
 Let $l=$ $\max_{i}\deg
h_i$, $V=$ ${\La}_{l}$,  $Y=$ ${\k} h_1+\ldots$
$+{\k} h_k\subset$ $V$   and $r=\dim
(I_{\mu})_{l}$. Note that by Lemma
\ref{lemma-regular-in-complete-intersection},
(\ref{item-poly-graded-does-not-depend}), $r$
does not depend on $\mu\in \Sp \Ga$ and hence
there is a well defined  map  $F: \Sp \Ga\to
Gr_{r}(\La_l)$ such that $F(\mu)=$
$(I_{\mu})_{l}$.  It follows from  Lemma
\ref{lemma-regular-in-complete-intersection},
(\ref{item-poly-min-present}) that $F$ is a
regular map.  Then by linear independency of
$h_1, \ldots h_k$ over $I$ it follows that $d=$
$\ds {\min_{\mu\in\Sp{\Ga}}}\dim Y\cap
(I_{\mu})_{l}=0$. Applying Lemma
\ref{lemma-min-intersection-is-open} we
conclude that  the set ${\sU}_1$ consisting of
$\mu$'s, such that $\dim Y\cap(I_{\mu})_{l}=0$,
is open dense in $\Sp{\Ga}$ proving the
statement
(\ref{item-linearly-indep-over-i-mu}). By Lemma
\ref{lemma-regular-over-ci-allmost-everywhere-reduced}
there exists an open dense set
${\sU}_{\sg}\subset \Sp \Ga$ such that
$I_{\mu}$ is radical for every  $\mu\in
{\sU}_{\sg}$. For such $\mu$'s  the linear
independence of $h_1,\ldots,h_k$ over $I_\mu$
is equivalent to the  linear independence of
the restrictions
 $h_1|_{{\sV}(I_{\mu})}$, $\ldots$, $h_k|_{{\sV}(I_{\mu})}$  as regular
functions. Thus the statement
(\ref{item-linearly-indep-as-functions})
follows from
(\ref{item-linearly-indep-over-i-mu})  and Lemma
\ref{lemma-regular-over-ci-allmost-everywhere-reduced}
if we set ${\sU}_2=$ ${\sU}_1\cap {\sU}_{\sg}$.
\end{proof}

\begin{lemma}{}\label{lemma-freedom-in-commutative-case}
Let $\La$ be a graded \CI algebra, $\sg=\{g_1,
\ldots, g_t\}$  a complete intersection of
homogeneous elements in $\La$, $I=(g_1,\ldots,
g_t)$ and $\Ga=\k [g_1, \ldots, g_t]$. Suppose
$I=$ $\ds \sum_{i=0}^{\infty}I_{(i)}$ is a
graded decomposition of $I$ and $H=$ $\ds
\sum_{i=0}^{\infty}H_{(i)}$ is a graded
complement of $I$ in $\La$ as a ${\k}-$vector
space. Then the mapping
$\pi:{\Ga}\otimes_{{\k}}H{\longrightarrow}$
${\La}$ defined by $\pi(\gamma\otimes h)=$
$\gamma h$, $\gamma\in{\Ga}$, $h\in H$ is an
isomorphism of $\Ga$-modules. In particular,
${\La}$ is a free module over ${\Ga}$.
\end{lemma}

\begin{proof} Note that $H_{(0)}={\k}$ and $\im\pi$ is a
${\Ga}-$submodule in ${\La}$ containing ${\Ga}$
and $H$. We prove by induction on $i$ that
${\La}_{(i)}=I_{(i)}+H_{(i)}\subset$ $\im \pi$.
If $f\in{\La}_{(i)}$ then $f=$ $f_I+f_H$,
$f_I\in I_{(i)}$, $f_H\in H_{(i)}\subset \im
\pi$. We have that $f_I=$ $\ds
{\sum_{j=1}^t}f_jg_j$, where $\deg f_j<i$,
$j=1,$ $\dots,$ $t$. By induction, $f_j\in$
$\im\pi$ for all $j$ and, since $\im \pi$ is a
$\Ga$-module, $f_I\in\im\pi$, therefore
$f=f_I+f_H\in\im \pi$. It is left to show that
$\pi$ is a monomorphism. Suppose
$h_1,\ldots,h_k\in$ $H$ are linearly
independent over ${\k}$, hence
 over
$I$, and for some
$\gamma_1,\ldots,\gamma_k\in$ ${\Ga}$ holds
$\gamma_1 h_1+\ldots$ $+\gamma_k h_k=0$.
 By
Lemma
\ref{lemma-equality-modulo-ideal-and-equality-as-functions},
(\ref{item-linearly-indep-as-functions}) there
exists an open dense set ${\sU}_2\subset$
$\Sp{\Ga}$ such that for each $\mu\in$
${\sU}_2$ the restrictions of functions
$h_1|_{{\sV}(I_{\mu})}$, $\ldots$,
$h_k|_{{\sV}(I_{\mu})}$ are linearly
independent over ${\k}$, where
$I_{\mu}=(g_1-g_1(\mu), \ldots, g_t-g_t(\mu))$.
Since $\gamma_i|_{{\sV}(I_{\mu})}=$
$\gamma_i(\mu)\in\k$, $i=1,\ldots,k$, we get
that $\gamma_1|_{{\sV}(I_{\mu})}=$ $\ldots$
$=\gamma_k|_{{\sV}(I_{\mu})}=0$ for every $\mu
\in {\sU}_2$. This implies
$\gamma_1|_{p_{\sg}^{-1}({\sU}_2)}=$ $\ldots$
$=\gamma_k|_{p_{\sg}^{-1}({\sU}_2)}=0$ where
$p_{\sg}$ is defined as in Lemma
\ref{lemma-regular-in-complete-intersection},
(\ref{item-poly-map-is-epi}). Since
$p_{\sg}^{-1}({\sU}_2)$ is dense in $\Sp\La$ we
conclude that $\gamma_1=$ $\ldots$
$=\gamma_k=0$.
\end{proof}

Now we are in the position to  prove our main result.

\bigskip

\noindent {\bf Proof of Theorem
\ref{theorem-main-analogue-kostant-theorem}.}
We  prove the statement for $U$ as a left
module (right module structure is treated
analogously). We apply Lemma
\ref{lemma-freedom-in-commutative-case} to
${\La}=$ ${\oU}$ and the sequence
$\{{\og}_1,\ldots,{\og}_t\}$.  Let $\oI$ be the
left ideal of $\oU$ generated by $\og_1,
\ldots, \og_t$, ${\oGa}=$
${\k}[{\og}_1,\ldots,{\og}_t]$ and ${\overline
H}=$ $\ds {\sum_{i=0}^{\infty}}$ ${\overline
H}_i$ be a graded complement to ${\oI}$ in
${{\La}}$. Then the map ${\overline \pi}:$
${\oGa}\otimes_{{\k}}{\overline
H}{\longrightarrow}$ ${\La}$ which sends
${\overline  \gamma} \otimes {\overline
h}\longmapsto$ ${\overline  \gamma} {\oh}$ is
an isomorphism of vector spaces. Choose for all
$i\geqs0$ the $\k$-vector spaces $H_i\subset$
$U_{i}$  such that ${\gr}:U{\longrightarrow}$
${\oU}$ induces a ${\k}-$linear isomorphism
$H_i$ onto ${\overline  H}_i$ and $H=$ $\ds
{\sum_{i=0}^{\infty}}H_i$. The filtrations on
$\Ga$ and $H$ induce a filtration on
${\Ga}\otimes_{{\k}}H$ by
$$({\Ga}\otimes_{{\k}}H)_k= \sum_{i=0}^n
{\Ga}_i\otimes_{{\k}}H_{k-i},\text{ and  in
particular } \gr ({\Ga}\otimes_{{\k}}H)_{(k)} =
\sum_{i=0}^n \oGa_i\otimes_{{\k}}\overline
H_{k-i}.$$ The map $\pi:$
${\Ga}\otimes_{{\k}}H{\longrightarrow}$ $U$,
  $\gamma \otimes h\longmapsto$
$\gamma h$ preserves the filtrations and the
induced graded map coincides with $\overline
\pi$. Since $\overline \pi$ is an isomorphism,
we see that $\pi$ is a $\k$-isomorphism  which
completes the proof.

\section{Applications}
\label{section-applications}

\subsection{Semisimple Lie algebras}
\label{subsectionj-semisimple-lie-algebras} Let
$U=\U(\g)$ be the universal enveloping algebra
of a semisimple Lie algebra of   rank $r$,
$I(\g)$
 the algebra of $\g$-invariants in $S(\g)$ (\cite{Di}, 2.4.11).
Then in $I(\g)$ there exists a system of
homogeneous algebraically independent
generators $f_1,\dots,f_r$ (\cite{Di}, Theorem
7.3.8) and the variety $\sV(f_1,\dots,f_r)$ is
irreducible of dimension $\dim \g- r$
(\cite{Di}, Theorem 8.1.3). If
$\varphi:S(\g)\to \U(\g)$ is the symmetrization
map  then it maps $I(\g)$ isomorphically onto
the center $Z(\g)$  of $\U(\g)$  and
$\sf=\{\varphi(f_1),\dots,\varphi(f_r)\}$ are
the generators of $Z(\g)$. Since
$\gr(\U(\g))\simeq S(\g)$ and
$\gr(\varphi(f_i))=f_i,$ $i=1,\dots,r$ the
classical Kostant theorem (\cite{Di}, 8.2.4)
follows from Theorem
\ref{theorem-main-analogue-kostant-theorem}
applied to $U=\U(\g)$ and the family $\sf$ of
generators of $Z(\g)$ and Proposition
\ref{proposition-bh},
(\ref{item-proposition-equidimensionality-is-equivalent-to-complete-intersection}).

\subsection{Restricted Yangians}
\label{subsectionj-restricted-yangians}

Let $p$ be a positive integer.
The {\it level $p$ Yangian\/} $Y_p(\gl_n)$ for the Lie algebra
$\gl_n$ (\cite{D}, \cite{C}) can be defined as the associative algebra
with generators $t_{ij}^{(1)}, \dots, t_{ij}^{(p)}$, $i,j=1,\ldots, n$,
subject to the relations
\begin{equation}
[T_{ij}(u), T_{kl}(v)]=\frac{1}{u-v}(T_{kj}(u)\hspace{1pt}
T_{il}(v)-T_{kj}(v)\hspace{1pt} T_{il}(u)),
\end{equation}
where $u,v$ are formal variables and
\begin{equation}
T_{ij}(u)=\delta_{ij}\, u^p+\sum_{k=1}^p t_{ij}^{(k)}\, u^{p-k}\in
\Y_p(\gl_n)[u].
\end{equation}
 These relations are equivalent to the following
\begin{equation}
[t^{(r)}_{ij}, t^{(s)}_{kl}] =\sum_{a=1}^{\min(r,s)}
\big(t^{(a-1)}_{kj} t^{(r+s-a)}_{il}-t^{(r+s-a)}_{kj} t^{(a-1)}_{il}\big),
\end{equation}
where $t^{(0)}_{ij}=\delta_{ij}$ and
$t_{ij}^{(r)}=0$ for $r\geqs p+1$.

The importance of the restricted Yangian
$Y_p(\gl_n)$ is motivated by the fact that any
irreducible finite-dimensional representation
of the full Yangian $Y(\gl_n)$ is a
representation of the restricted Yangian for
some $p$ \cite{D2}. Note that the level $1$
Yangian  $\Y_1(\gl_n)$ coincides with the
universal enveloping algebra $\U(\gl_n)$.

Set $\deg t_{ij}^{(k)}=k$. This defines a
filtration on $\Y_p(\gl_n)$. The following
analogue of the Poincar\'e--Birkhoff--Witt
theorem for the algebra $\Y_p(\gl_n)$
(\cite{C}, \cite{Mo}), shows that $\Y_p(\gl_n)$
is a special filtered algebra.

\begin{proposition}\label{thm:pbw}
The associated graded algebra $\overline {\Y}_p(\gl_n)=\gr \Y_p(\gl_n)$ is
a polynomial algebra in variables $\overline {t}_{ij}^{(k)}$, $i,j=1, \ldots,
n$, $k=1, \ldots, p$.
\end{proposition}

Note that the  grading on $\overline \Y_p(\gl_n)$ induced by the
grading on
 $\Y_p(\gl_n)$,
 $\deg \overline {t}_{ij}^{(k)}=k$, does not coincide with the standard
 polynomial grading for $p>1$. Set
$\T(u)=(T_{ij}(u))_{i,j=1}^n$
and consider the following element in $\Y_p(\gl_n)[u]$, called
 \emph{quantum determinant}
\begin{equation}
\begin{split}
\qdet \T(u)= \sum_{\sigma\in S_n} \sgn
(\sigma)\,
T_{1{\sigma(1)}}(u)T_{2{\sigma(2)}}(u-1)...T_{n{\sigma(n)}}(u-n+1).  \\
\end{split}
\end{equation}

The  coefficients $d_{s}$ by the powers
$u^{np-s}$, $s=1, \ldots, np\,$  of  $\,\qdet
\T(u)$ are algebraically independent generators
of the center  of $\Y_p(\gl_n)$ (\cite{C},
\cite{Mo}).

For $F=\sum_i f_iu^i\in \Y_p(\gl_n)[u]$ denote
$\overline F=\sum_i \overline {f}_i u^i\in
\overline {\Y}_p(\gl_n)[u]$. Also we denote
$X_{ij}^{(k)}=\overline {t}_{ij}^{(k)}$,
$X_{ij}(u)=\overline {T}_{ij}(u)$ and
$X(u)=(X_{ij}(u))_{i,j=1}^n$. Since $\overline{
{T}_{ij}(u-\lambda)}=X_{ij}(u)$ for any
$\lambda\in\k$, one can easily check  that $\gr
\qdet T(u)=\det X(u)$.

\begin{lemma}\label{lemma-center-of-restricted-yangian-gives-regular-sequence}
The sequence $\overline d_1, \ldots, \overline d_{np}$ is a
complete intersection for $\overline \Y_p(\gl_n)$.
\end{lemma}

\begin{proof}
Due to Proposition \ref{proposition-bh},
(\ref{item-subsequeqnce-is-complete-intersection})
it is enough to prove that the sequence
$$X_{ij}^{(k)}, i\neq j, i,j=1, \ldots, n, k=1, \ldots p; \ \overline d_1,
\ldots, \overline d_{np}$$ is complete intersection for $\overline
\Y_p(\gl_n)$.
 Let $c_s$ is a polynomial in variables $X_{11}^{(1)}$, $X_{11}^{(2)}$, $\ldots$,
 $X_{11}^{(p)}$, $X_{22}^{(1)}$, $\ldots$,
$X_{nn}^{(p)}$ obtained from $\overline d_s$ by
substituting $X_{ij}^{(k)}=0$ for all $i\neq j$
and all $k$. Due to Lemma
\ref{lemma-substituion-0-in-regular-gives-regular}
we only need to show the regularity of the
sequence $c_1, \ldots, c_{np}$ in the
polynomial ring $\k[X_{ii}^{(k)}, i=1,
\ldots,n, k=1,\ldots, p],$ hence by Proposition
\ref{proposition-bh},
(\ref{item-proposition-equidimensionality-is-equivalent-to-complete-intersection}),
we should show, that the variety $Z=\sV(c_1,
\ldots, c_{np})$
 is equidimensional of dimension $0$.

Substituting $X_{ij}^{(k)}=0$ into the matrix
$X(u)$ we obtain  that $c_s$ is the coefficient
by $u^{np-s}$ in $\det \diag\{X_{11}(u),
\ldots, X_{nn}(u)\}=\prod_{i=1}^n X_{ii}(u)$.
Consider the regular map $\varphi: \k^{np}\to
\k^{np}$ which sends $(a_{ii}^{(k)})\in
\k^{np}$ to the coefficients of the following
monic polynomial
$\prod_{i=1}^n(u^p+a_{ii}^{(1)}u^{p-1}+ \ldots
+ a_{ii}^{(p)})$. Obviously, $Z=
\varphi^{-1}(0)$.  Since $\k[u]$ is a factorial
domain, $u^{np}=\underbrace{u^p\cdot u^p\cdot
\ldots\cdot u^p}_{n}\ $ is a unique
decomposition of $u^{np}$ in the products of
monic polynomials of degree $p$ and  hence
$Z=\{0\}$. This completes the proof.
\end{proof}

Applying Theorem \ref{theorem-main-analogue-kostant-theorem} and
Lemma
\ref{lemma-center-of-restricted-yangian-gives-regular-sequence} we
 obtain the following analogue of the Kostant theorem
for the restricted Yangians.

\begin{theorem}
\label{theorem-application-to-restricted-yangians}
For all $n,p\geqs 1$ the restricted Yangian
$\Y_p(\gl_n)$ is a free module over its center.
\end{theorem}

\subsection{Current algebras}
\label{subsection-current-algebra}

In this section we consider the polynomial
current Lie algebra $\g\,
 =\gl_n(\C)\otimes \C[x]$  and its
restricted quotient algebra $\g_m=\g_m(n)$,
$m>0$, by the ideal $\sum_{k\geqs m}\gl_n
\otimes x^k$ (\cite{RT}). We will show that the
universal enveloping algebra $\U(\g_m)$ is a
free module over its center for any $m>0$. This
generalizes in the case of $\gl_n$ the result
of Geoffriau for  Takiff algebra $\g_2$
(\cite{G1},\cite{G2}).

In \cite{Mo} were constructed the families of
algebraically independent generators of the
center of $\U(\g_m)$. Following \cite{Mo} let
$E_{ij}$, $i,j=1, \ldots, n$ be the standard
basis of
 $\gl_n$, $E_{ij}^{(k)}=E_{ij}\otimes x^k$ with $1\leqs i,j\leqs n$,
$0\leqs k\leqs m-1$   a basis of $\g_m$. Set

$$F_{ij}^{(r)}=E_{ij}^{(r-1)},\,1<r\leqs m \text{
and }
F_{ij}^{(1)}=E_{ij}^{(0)}-m(j-1)\delta_{ij}.$$

For each $k\in \{1, \ldots,  mn\}$ let $r\in \{1, \ldots, m\}$ and $s\in \{1, \ldots, n\}$ be such that
$k=m(s-1)+r$.  Then the elements

\begin{equation}
\xi_k= \sum_{\substack{i_1<\ldots <i_s\\ j_1+ \ldots +j_s=k}}\sum_{\sigma\in
S_s}\sgn(\sigma)F_{i_{\sigma(1)}i_1}^{(j_1)}\ldots F_{i_{\sigma(s)}i_s}^{(j_s)},
\end{equation}
are algebraically independent generators of the
center  of $\U(\g_m)$ (\cite{Mo}). Note that
$\U(\g_m)$ is a special filtered algebra with
respect to the standard grading, $\deg
E_{ij}^{(k)}=1$ for all $i,j,k$.

We will show that the sequence $\overline \xi_1, \ldots, \overline
\xi_{mn}$ is complete intersection for  $\oU(\g_m)$ where

\begin{equation}\label{tsentr cur}
\overline \xi_k= \sum_{\substack{i_1<\ldots <i_s\\ j_1+ \ldots +j_s=k}}\sum_{\sigma\in
S_s}\sgn(\sigma)\overline F_{i_{\sigma(1)}i_1}^{(j_1)}\ldots \overline F_{i_{\sigma(s)}i_s}^{(j_s)}.
\end{equation}
As in the case of the Yangian $\Y_p(\gl_n)$ we
complete this sequence by the elements
$\overline F_{ij}^{(l)}$ for all $i,j=1,
\ldots, n$, $i\neq j$ and $l=1, \ldots, m$ and
apply Lemma
\ref{lemma-substituion-0-in-regular-gives-regular} by substituting
 $\overline F_{ij}^{(l)}=0$.
Hence it is enough to prove the regularity of
the sequence $\gamma_1^m, \ldots,
\gamma_{mn}^m$, where
 $\ds \gamma_k^m=\sum_{\substack{i_1<\ldots <i_s\\ j_1+ \ldots +j_s=k}}
\overline F_{i_{1}i_1}^{(j_1)}\ldots \overline
F_{i_si_s}^{(j_s)}$ (compare with (3.2) in
\cite{Mo}). We will show that
$\sX=\sV(\gamma_1^m, \ldots,
\gamma_{mn}^m)=\{0\}$ by induction on $m$.
Suppose that $m=1$. Since for each $k=1,
\ldots, n,$ $\gamma_k^1$ is the elementary
symmetric polynomial of degree $k$ in variables
$\overline F_{11}^{(1)}, \ldots, \overline
F_{nn}^{(1)}$ we have that $\sV(\gamma_1^1,
\ldots, \gamma_{n}^1)=\{0\}$. Suppose now that
the sequence $\gamma_1^{m-1}, \ldots,
\gamma_{(m-1)n}^{m-1}$ is a complete
intersection for $\oU(\g_m)$ and
$\sV(\gamma_1^{m-1}, \ldots,
\gamma_{(m-1)n}^{m-1})=\{0\}$. Note that
$\gamma_{ms}^m$ is the elementary symmetric
polynomial of degree $s$ in variables
$\overline F_{11}^{(m)}, \ldots, \overline
F_{nn}^{(m)}$, $s=1, \ldots, n$. Hence
$\overline F_{11}^{(m)}(x)=\ldots=\overline
F_{nn}^{(m)}(x)=0$ for every $x\in \sX$.
Substituting these values in $\gamma_k^m$ for
every $k$ that does not divide $m$ we obtain
the sequence $\gamma_1^{m-1}, \ldots,
\gamma_{(m-1)n}^{m-1}$ which is a complete
intersection by our assumption and therefore
$\sX=0$. By induction on $m$ we conclude that
the sequence $\overline \xi_1, \ldots,
\overline \xi_{mn}$ is a complete intersection
for $\oU(\g_m)$.

 Applying Theorem
\ref{theorem-main-analogue-kostant-theorem} we
immediately obtain the following analogue of
the Kostant theorem for restricted current
algebras:

\begin{theorem}
\label{theorem-application-to-current-algebras}
For all $m,n\geqs 1$ the algebra $\U(\g_m(n))$
is a free module over its center.
\end{theorem}

\section{Acknowledgment}
The first author is a Regular Associate of the ICTP and is
supported by the CNPq grant (Processo 300679/97-1). He is
grateful to the University of Sydney for support and hospitality
during his visit.
  The second
author is grateful to FAPESP
 for the financial support (Processo 2002/01866-7)  and
to the
University of S\~ao Paulo for  hospitality
during his visit. Both authors are deeply grateful to Alexander Molev for
posing this problem for restricted Yangians and current algebras,
his encouragement and many useful discussions.


\begin{thebibliography}{20}


\bibitem[B]{B}{Bourbaki, N.}, \textit{Alg\'ebre homologique, Alg\'ebre chap. X},
Masson, 1980.


\bibitem[C]{C}
{Cherednik I.V.},
\textit{Quantum groups as hidden symmetries of classic
representation theory},
In: ``Differential Geometric Methods in Physics"
(A. I. Solomon, Ed.),
World Scientific,
Singapore,
1989,
 47--54.

\bibitem[Di]{Di}{Dixmier J.}, \textit{Algebras enveloppantes},
Gauthier-Villarrd, Paris, 1974.


\bibitem[D1]{D}
{Drinfeld V.G.},
\textit{Hopf algebras and the
quantum Yang--Baxter equation},
{Soviet Math. Dokl.} \textbf
{32} (1985), 254--258.


\bibitem[D2]{D2}
{Drinfeld V.G.},
\textit{A new realization of Yangians and quantized affine algebras},
{Soviet Math. Dokl.} \textbf
{36} (1988), 212--216.

\bibitem[DFO]{DFO}{Drozd Yu.A., Ovsienko S.A., Futorny
V.M.}, \textit{Harish - Chandra subalgebras and
Gelfand Zetlin modules,} in: "Finite
dimensional algebras and related topics", NATO
ASI Ser. C., Math. and Phys. Sci.,
\textbf{424}, (1994), 79-93.



\bibitem[Ei]{Ei}{Eisenbud D.}, \textit{Commutative algebra}, Graduate
Texts in Math., \textbf{150}, Springer, 1994.


\bibitem[G1]{G1}
{Geoffriau F.},
\textit{Une propri\'et\'e des alg\`ebres de Takiff},
C. R. Acad. Sci. Paris \textbf{319} (1994), S\'erie I, 11--14.

\bibitem[G2]{G2}
{Geoffriau F.},
\textit{On the center of the enveloping algebra of a Takiff algebra},
Ann. Math. Blaise Pascal, \textbf{1} (1995), 15--31.

\bibitem[K]{K}{Kostant B.}, \textit{Lie groups representations on polynomial
rings}.
Amer.J.Math. \textbf{85}, (1963), 327-404.

\bibitem[Ma]{Ma} {Matsumura H.}, \textit{Commutative ring theory},
Cambridge Studies in Advanced Mathematics, 8, Cambridge University
Press, 1997.

\bibitem[Ma1]{Ma1} {Matsumura H.}, \textit{Commutative Algebra},
Math. Lecture Note Series, W.A.Benjamin, Inc., New York,
 1970.

\bibitem[Mo]{Mo}
{Molev A.I.}, \textit{Casimir elements for certain polynomial
current Lie algebras}, In: ``Group 21, Physical Applications and
Mathematical Aspects of Geometry, Groups, and Algebras," Vol. 1,
(H.-D.~Doebner, W.~Scherer, P.~Nattermann, Eds). World Scientific,
Singapore, 1997, 172--176.


\bibitem[MNO]{MNO}
{Molev A., Nazarov M, Olshanski G.},
\textit{Yangians and classical Lie algebras}, Russian Math. Surveys
\textbf{51}:2 (1996), 205-282.


\bibitem[O1]{O1}{Ovsienko S.}, \textit{Strongly nilpotent matrices and
Gelfand-Tsetlin modules}, J.Linear Algebra and
Appl., (2003), 19pp. to appear.

\bibitem[O2]{O2}{Ovsienko S.}, \textit{Finiteness statements for
Gelfand-Tsetlin modules}, In: Proceedings of International Algebraic
Conference, Sumi, 2002.

\bibitem[RT]{RT}{Rais M., Tauvel P.}, \textit{Indice et polynomes invariants pour
certaines algebres de Lie}, J. Reine Agnew. Math. \textbf{425} (1992), 123-140.

\bibitem[W]{W} {Ch. A. Weibel},
\textit{ An introduction to homological
algebra}, Cambridge Studies in Advanced
Mathematics {\bf 38}, Cambridge University
Press (1994).

\end{thebibliography}
\end{document}